\documentclass{amsart}
\usepackage{amssymb}
\usepackage{amscd}
\usepackage{epsfig}
\newtheorem{theorem}{Theorem}

\newtheorem{proposition}[theorem]{Proposition}

\theoremstyle{definition}
\newtheorem{definition}{Definition}

\numberwithin{equation}{section}

\newcommand{\bc}{\mathbb{C}}
\newcommand{\bp}{\mathbb{P}}
\newcommand{\bz}{\mathbb{Z}}

\newcommand{\bh}{\mathbb{H}}
\newcommand{\Z}{\mathbf{Z}}

\newcommand{\cs}{\mathcal{S}}

\newcommand{\hk}{\hookrightarrow}

\newcommand{\med}{\medskip}
\newcommand{\la}{\longrightarrow}
\newcommand{\bfl}{\begin{flushleft}}
\newcommand{\efl}{\end{flushleft}}

\newcommand{\ltm}{LM^{-TM}}
\newcommand{\mtm}{M^{-TM}}
\newcommand{\bcp}{\bc \bp}

\newcommand{\lmm}{LM \times_M LM}
\newcommand{\evi}{\text{ev}_\infty}
\newcommand{\evl}{ev}

\newcommand{\ttt}{\tilde{t}}

\newcommand{\into}{\hookrightarrow}
\title{The Loop Homology Algebra of Spheres and Projective Spaces}

\author[R.L. Cohen]{Ralph L. Cohen}
\address{Dept. of Mathematics\\
Stanford University\\
Stanford, California 94305}
\email[Cohen]{ralph@math.stanford.edu}
\thanks{The first author was partially supported by a grant from the NSF} 

\author[J.D.S Jones]{John D.S. Jones}
\address{Department of Mathematics\\
University of Warwick\\
Coventry, CV4 7AL England}
\email[Jones]{jdsj@maths.warwick.ac.uk}

\author[J. Yan]{Jun Yan}
\address{Dept. of Mathematics\\
Stanford University\\
Stanford, California 94305}
\email[Yan]{jyan@math.stanford.edu}

\date{\today}

\begin{document}

\begin{abstract}
  In \cite{CS} Chas and Sullivan defined an intersection product on the homology $H_*(LM)$ of the 
space of smooth loops in a closed, oriented manifold $M$. In this paper we will  use the homotopy theoretic realization of this product 
described by the first two authors in \cite{CJ} to construct a  second quadrant spectral sequence of
algebras converging to the loop homology multiplicatively, when $M$ is simply connected.  The $E_2$ term of this spectral sequence is $H^*(M;
  H_*(\Omega M))$ where the product is given by the   cup product on the cohomology of the manifold
$H^* (M)$ with coefficients in the Pontryagin ring structure on the homology of its based loop space $H_*(\Omega M)$.  We then use this spectral
sequence to compute the ring structures of $H_* (LS^n)$ and $H_* (L\bcp^n)$.  
\end{abstract}
\maketitle

\section*{Introduction}
The loop homology of a closed orientable manifold $M^d$ of degree $d$ is the ordinary homology of the free loop space $LM = \text{Map}(S^1, M^d)$, with degree shifted by $-d$, i.e.
$$\bh_*(LM; \bz) = H_{*+d}(LM; \bz).$$
In \cite{CS}, Chas and Sullivan defined a type of intersection product   on the chains of
$LM$, yielding an algebra structure on $\bh_*(LM)$.  

Roughly, the loop product is defined
  as follows.   Let
$\alpha : \Delta ^p \to LM$ and $\beta : \Delta^q \to LM$ be  singular 
simplices in $LM$.  The evaluation at $1 \in S^1 \subset \bc$ defines a map $ev :LM \to M$.  Assume that $ev \circ \alpha$ and $ev \circ
\beta $ define a map $$\Delta^p \times \Delta ^q \to LM \times LM \to M \times M$$ that is transverse to the diagonal.  At each point $(s,t)
\in \Delta^p \times \Delta^q$  where $ev \circ \alpha$ intersects $ev \circ \beta$, one can define a single loop by first traversing the loop
$\alpha (s)$ and then traversing the loop $\beta (t)$.  This then defines a chain $\alpha \circ \beta \in C_{p+q-d}(LM)$.   In \cite{CS} Chas
and Sullivan showed  that this procedure defines a chain map
$$
C_p(LM) \otimes C_q(LM) \to C_{p+q-d}(LM)
$$
which induces an associative, commutative algebra structure on  the loop homology, $\bh_*(LM)$. 
Chas and Sullivan also described other structures this pairing induces, such as a Lie algebra
structure on the equivariant homology of the loop space.   In
\cite{CJ}, the first two authors used the Pontryagin - Thom construction to show that the loop
product is realized on the homotopy level on Thom spaces (spectra)  of bundles over the loop
space.   In particular let  $TM$ denote the tangent bundle of $M$, and $-TM$ denotes its inverse
as a  virtual bundle in $K$ - theory.  Let $\mtm$ denote the Thom spectrum of this bundle, and
$\ltm$ the Thom spectrum of
$ev^*(-TM)$.  Then in \cite{CJ} it was shown that $\ltm $ is a homotopy commutative ring spectrum with unit, whose product realizes the Chas -
Sullivan product in homology, after applying the Thom isomorphism, $H_q(\ltm) \cong H_{q+d}(LM) \cong \bh_{q}(LM)$. 

The goal of this paper is to describe a spectral sequence of algebras converging to the loop homology algebra of  a manifold, and to use it to
compute the loop homology algebra of spheres and projective spaces.  More specifically we shall prove the following theorems.
\med
\begin{theorem}\label{ss}
Let $M$ be a closed, oriented, simply connected  manifold.  There is a second quadrant spectral sequence of
algebras
$\{E_{p,q}^r, d^r \, : p \leq 0, q \geq 0\}$ such that
\begin{enumerate}
\item 
$E^r_{*,*}$ is an algebra and the
differential
$d^r : E^r_{*,*}
\to E^r_{*-r, *+r-1}$ is a derivation for each $r\geq 1$. 
\item The spectral sequence converges  to the loop homology $\bh_*(LM)$ as algebras.  That is,
$E^\infty_{*,*}$ is the associated graded algebra to a natural filtration of the algebra $\bh_*(LM)$. 
 \item  For $m,n \geq 0$,   $$E^2_{-m, n } \cong H^m(M; H_n(\Omega M)).$$ Here $\Omega M$ is the space of base point
preserving loops in $M$.  Furthermore the isomorphism
$E^2_{-*,*} \cong H^*(M; H_*(\Omega M))$ is an isomorphism of algebras, where the algebra structure on 
$ H^*(M; H_*(\Omega M))$ is given by the cup product on the cohomology  of $M$ with coefficients in the Pontrjagin ring
$H_*(\Omega M)$.    
 \item The spectral sequence is natural with respect to smooth maps between manifolds.
\end{enumerate}
\end{theorem}
\med
\noindent We then use this spectral sequence to do the following calculations. Let $\Lambda [x_1,  \cdots , x_n]$ denote the
exterior algebra (over the integers)  generated by $x_1, \cdots x_n$, and let $\bz[a_1, \cdots , a_m]$ denote the polynomial algebra  generated
by $a_1, \cdots, a_m$.  

\med
\begin{theorem}\label{spheres} There exist isomorphisms of graded algebras,
\begin{enumerate}
\item $$\bh_{*}(LS^1) \cong \Lambda [a] \otimes \bz[t, t^{-1}]$$
where $a \in \bh_{-1}(LS^1)$   and $t,  t^{-1} \in \bh_0 (LS^1)$.   
\item For $n >1$, 
$$\bh_*(LS^n) = \begin{cases}
\Lambda[a] \otimes   \bz[u], &\text{for $n$ odd}\\
(\Lambda[b] \otimes \bz[a, v]) / (a^2, ab, 2av), &\text{for $n$ even},\end{cases}
$$
where $a \in \bh_{-n}(LS^n)$, $b \in \bh_{-1}(LS^n)$, $u \in \bh_{n-1}(LS^n)$,  and $v \in \bh_{2n-2}(LS^n)$.  
\end{enumerate}
\end{theorem}

\med 
\begin{theorem}\label{projspace}
There is an isomorphism  of algebras,
$$\bh_* (L\bcp^n) \cong (\Lambda[w] \otimes \bz[c, u]) /(c^{n+1}, (n+1)c^n u, wc^n)$$
where $w \in \bh_{-1}(L\bcp^n)$, $c \in \bh_{-2}(L\bcp^n)$, and $u \in \bh_{2n}(L\bcp^n)$.
\end{theorem}

\med
\noindent The organization of this paper is as follows. 
In section 1 we will review the construction of the loop product  that was used in \cite{CJ} and
describe it on the  chain level. In section 2 we will construct the spectral sequence and prove
theorem \ref{ss}. In section 3 we use this spectral sequence to do the calculations presented in
theorems
\ref{spheres} and \ref{projspace}. 
\med
\section{The Loop Product}

For the remainder of the paper, let $M$ be a closed, connected,  oriented manifold of dimension
$d$. The goal of this section is to give a description of the loop product on $\bh_* (LM)$  at the
chain level.  Of course this was done originally by Chas and Sullivan in \cite{CS}, but our approach
will be slightly different.  It will be more amenable to the construction and the analysis of the
loop homology spectral sequence to be done in the next section.

We begin by recalling a
chain level description of  the intersection product arising from Poincare duality on the homology of
the manifold,
$$
<,> : H_q(M) \otimes  H_p(M) \to H_{p+q-d}(M).
$$

For  pairs of spaces
$(X,A) $ with $A
\subset X$, denote by   $C_* (X, A)$ and  $C^* (X, A)$ the groups of singular
chains and cochains.   

Recall that the normal bundle of the diagonal embedding $\Delta : M \to M \times M $
  is naturally isomorphic to the tangent
bundle
$TM$.  Let
$M^{TM}$ denote the Thom space of  this bundle.  The   Thom-Pontryagin map for the diagonal
embedding is therefore a map
\begin{equation}\label{tau}
\tau: M \times M \to M^{TM},
\end{equation}
which collapses everything outside  a tubular neighbourhood of the diagonal to the base point of
$M^{TM}$.

Choose a Riemannian metric on $TM$ and let $D_TM$ and $S_TM$ be the unit disk and sphere
bundles of $TM$. Since $S_TM$ has a   collar,
  $S_TM \times I \into D_TM$, there is a chain equivalence
$$
\theta_\sharp: C_*(D_TM/S_TM, *) \to C_*(D_TM \cup cS_TM, cS_TM) \to C_*(D_TM, S_TM)
$$
  where $*$ is the base point of $D_TM/S_TM$,   $cS_TM$ is the cone on $S_TM$, and $D_TM \cup cS_TM$ is the mapping cone of the inclusion $S_TM
\into D_TM$.  Let
 $t \in C^d(D_TM, S_TM)$  be a cochain that represents the Thom class $[t] \in H^d (D_TM, S_TM)$. Then taking the  
cap product  with $t$ at the
chain level followed by the projection from
$D_TM$ to $M$,
\begin{equation}\label{sigma}
\begin{CD}
\sigma_\sharp: C_* (D_TM, S_TM) @>\cap t>> C_{*-d} (D_TM) @>\pi_*>> C_{*-d}(M)
\end{CD}
\end{equation}
induces the Thom isomorphism 
$$
\begin{CD}
{\tilde H}_* (M^{TM}) @>\theta_* >> H_*(D_TM, S_TM) @>\sigma_* >> H_{*-d} (M)
\end{CD}
$$
in homology.

Let $\epsilon_\sharp: C_*(D_TM/S_TM) \to C_*(D_TM/S_TM, *)$ be the projection  map onto the quotient
complex, and   denote by $u_\sharp$ the composition of the chain maps
\begin{equation}\label{usharp}
\begin{CD}
u_\sharp: C_*(M^{TM}) @>\epsilon_\sharp>> C_*(D_TM/S_TM, *) @>\theta_\sharp>> C_*(D_TM, S_TM) @>\sigma_\sharp>> C_{*-d}(M).
\end{CD}
\end{equation}

Now let $A$ be a graded ring.  Recall that the intersection product structure on $H_*(M; A)$ is
defined so that the Poincare duality isomorphism
\begin{align}
D: H^*(M;A) &\to H_{d-*}(M;A) \notag \\
\alpha &\to \alpha \cap [M] \notag
\end{align}
is an isomorphism of graded algebras.  Here $[M] \in H_d(M; A)$ is the fundamental class determined by the orientation of $M$.
The commutativity of the following
diagram  describes the well known relation between   the intersection product, the Thom - Pontryagin map and the Thom isomorphism:
$$
\begin{CD}
H_{p+q}(M \times M) @>\tau_*>> H_{p+q}(M^{TM}) @>u_*>> H_{p+q-d}(M) \\
@A\times AA & & @| \\
H_p (M) \otimes H_q (M) & @>(-1)^{d(d-p)}{<\cdot>}>> & H_{p+q-d}(M),
\end{CD}
$$
where $<\cdot>$ is the intersection product on $H_*(M)$,  and $\times $ denotes the cross product.   We therefore have the following  chain
description of the intersection pairing: 

\begin{proposition}\label{intersection}  Let $M$ be as above.  Then the following composition of chain maps
$$
\begin{CD}
C_p(M) \otimes C_q(M) @>\times >> C_{p+q}(M \times M) @> \tau_\sharp >> C_{p+q}(M^{TM}) @> u_\sharp >> C_{p+q-d}(M)
\end{CD}
$$ is a representative of  the intersection product
$$
(-1)^{d(d-p)}{<\cdot>}: H_p (M)  \otimes H_q (M) \to H_{p+q-d}(M).
$$
 \end{proposition}

\med
\noindent We next turn our attention to the loop space and the loop product.   Let $LM = C^\infty(S^1, M)$ be
the free loop space of
$M$, where $S^1$ denotes   the unit circle in the complex line, parametrized by
$exp: [0, 1] \to S^1$ with $1\in S^1$ chosen as the basepoint. We now recall some
constructions in \cite{CJ}.  

The loop evaluation map
\begin{align}
\evl: LM &\to M \notag\\ 
\gamma &\to \gamma (1), \notag
\end{align}
is a Serre fibration.    
Let  $\lmm$ denote the pull back of the product of this fibration with itself along the diagonal
embedding
$\Delta: M
\to M
\times M$, denoted by $\lmm$:  
$$
\begin{CD}
\lmm  @> \tilde \Delta>> LM \times LM \\
@V \evi VV  @VV \evl \times \evl V \\
M @>>\Delta > M \times M.
\end{CD}
$$
Notice that  $\lmm$ is the
space of pairs of loops having the same basepoint. 
The map $\evi$ in this fiber square is given by  $\evi(\alpha, \beta) = \alpha(1) = \beta(1)$.   The map $\tilde \Delta$ is an embedding
of a codimension $d$ infinite dimensional submanifold of the infinite dimensional manifold $LM \times LM$.

As shown in \cite{CJ} this pullback square allows for a Thom - Pontryagin map $\tilde \tau : LM \times LM \to  (\lmm)^{TM} $,
where $ (\lmm)^{TM} $ denotes the Thom space of the pull back bundle $\evi^*(TM)$, which is the normal bundle of  the embedding $\tilde \Delta$. 
Notice that we have a commutative diagram of Thom - Pontryagin maps:
$$
\begin{CD}
 LM \times LM@>\tilde \tau>> (\lmm)^{TM} \\
@V\evl \times \evl VV @V\evi VV \\
M \times M @>\tau>> M^{TM}.
\end{CD}
$$
In this diagram the map $\evi$ is  the induced map on the Thom spaces.    

Loop composition $\alpha \beta$ is defined for two loops $\alpha$ and $\beta$ having the same base point by first traversing the loop $\alpha$, then the loop $\beta$, i.e.
\begin{equation}\label{gamma}
\alpha \beta (t) = \begin{cases} \alpha(2t) \quad &\text{if} \ 0\leq t \leq \frac{1}{2}, \\
\beta(2t-1) \quad &\text{if} \ \frac{1}{2} \leq t \leq 1.\\
\end{cases}
\end{equation}
Denote this operation by
\begin{align}
\gamma : \lmm &\to LM \notag \\
(\alpha, \beta) & \to \alpha \beta. \notag
\end{align}

Notice that $\gamma$ preserves the base points of loops in $\lmm$ and $LM$, thus the composition
$$
\begin{CD}
\lmm @>\gamma>> LM @>\evl>> M
\end{CD}
$$
coincides with $\evi$.  Therefore $\gamma$ induces a map of bundles, $\gamma : \evi^*(TM) \to ev^*(TM)$, and 
therefore an induced map of
Thom spaces,  
$$\tilde \gamma: (\lmm)^{TM} \to LM^{TM}.$$
Putting the above maps together, we obtain the following commutative diagram
\begin{equation}\label{fibrations}
\begin{CD}
LM \times LM @>\tilde\tau>> (\lmm)^{TM} @>\tilde\gamma>> LM^{TM} \\
@V\evl \times \evl VV @V\evi VV @V\evl VV \\
M \times M @>\tau>> M^{TM} @>= >> M^{TM}.
\end{CD}
\end{equation}

Let $D_LM$ and $S_LM$ be the pulls back of $D_TM$ and $S_TM$ via $\evl$.  They are the unit   disk and sphere bundles of
$\evl^*(TM)$. Recall that the cochain $t$ in $C^d (D_TM, S_TM)$ represents the Thom 
class and therefore its pull back $\tilde t = \evl^*(t)$ is a cochain in $C^d (D_LM, S_LM)$ 
representing the Thom class of $\evl^*(TM)$. Similar to (\ref{sigma}),  capping with $\tilde t$ at
the chain level followed by the projection of
$D_LM$ to
$LM$
$$
\begin{CD}
\tilde\sigma_\sharp: C_* (D_LM, S_LM) @>\cap \tilde t>> C_{*-d} (D_LM) @>\tilde\pi_*>> C_{*-d}(LM)
\end{CD}
$$
induces the Thom isomorphism
$$
\begin{CD}
{\tilde H}_* (LM^{TM}) @> \cong >> H_*(D_LM, S_LM) @>{\tilde\sigma}_* >> H_{*-d} (LM)
\end{CD}
$$
in homology. Similar to the argument used for (\ref{usharp}), there is a chain map
\begin{equation}\label{ust}
\begin{CD}
\tilde{u}_\sharp: C_*(LM^{TM}) @>\tilde{\epsilon}_\sharp>> C_*(D_LM/S_LM, *) @>\tilde{\theta}_\sharp>> C_*(D_LM, S_LM) @>\tilde{\sigma}_\sharp>> C_{*-d}(LM),
\end{CD}
\end{equation}
such that the following diagram commutes:
$$
\begin{CD}
H_* (LM^{TM}) @> {\tilde u}_* >> H_{*-d} (LM) \\
@V \evl_* VV @V \evl_* VV \\
H_* (M^{TM}) @> u_* >> H_{*-d}(M).
\end{CD}
$$

In \cite{CJ}  the first two authors proved that the composition of $\tilde \tau_*$, $\tilde \gamma_*$, and $\tilde u_*$ realizes the Chas - Sullivan  loop
product.  That is, the following diagram commutes:  
\begin{equation}\label{prod}
\begin{CD}
H_p(LM) \otimes H_q(LM) & @>{\qquad\quad (-1)^{d(d-p)} (\circ)  \quad\qquad}>> & H_{p+q-d}(LM)  \\
@V \times VV & & @A\tilde u_* AA  \\
H_{p+q}(LM \times LM) @>\tilde \tau_* >> H_{p+q}(\lmm)^{TM} @>\tilde \gamma_* >> H_{p+q}(LM^{TM})
\end{CD}
\end{equation} where $\circ :  H_p (LM) \otimes H_q(LM) \to H_{p+q-d}(LM)$ is the Chas - Sullivan loop product. 
We therefore have the following proposition, which should be viewed as the analogue of   proposition \ref{intersection}. 

\begin{proposition}\label{loopprod} The composition of the four chain maps  
\begin{align}
\times:\quad & C_*(LM) \otimes C_*(LM) \to C_*(LM \times LM), \label{c1} \\
\tilde{\tau}_\sharp:\quad & C_*(LM \times LM) \to C_*(\lmm)^{TM}, \label{c2} \\
\tilde{\gamma}_\sharp:\quad & C_*(\lmm)^{TM} \to C_*(LM^{TM}), \text{\qquad and} \label{c3} \\
\tilde{u}_\sharp:\quad & C_*(LM^{TM}) \to C_{*-d}(LM). \label{c4} 
\end{align} gives a chain representative of the Chas - Sullivan loop product up to sign
\begin{align}
 H_p (LM) \otimes H_q(LM)  &\to H_{p+q-d}(LM) \notag \\
\alpha \otimes \beta &\to \alpha \circ \beta. \notag
\end{align}
\end{proposition}
\med 
\section{The Loop Algebra Spectral Sequence}

In this section we will describe the loop algebra spectral sequence and prove theorem \ref{ss}.  To do this we 
describe a filtration of simplicial sets arising from the fibration $\Omega M
\into LM \to M$, which will induce the Serre spectral sequence for this fibration.  We then analyze
how the loop product behaves with respect to this spectral sequence, using the chain level
description of the product given in the last section.  We then apply Poincare duality to regrade the
spectral sequence (and in particular change a first quadrant spectral sequence into a second
quadrant one), and prove theorem \ref{ss}. 

Given a topological space $X$, let $\cs. X$ denote its singular simplicial set.  The $p$ - simplices are given by
singular simplices $S_p(X) = \{\sigma : \Delta^p \to X\}$, where $\Delta^p$ is the standard $p$ - simplex,
$\Delta^p = \{(x_0, \ldots, x_p) \in \mathbb{R}^{p+1} : x_i \ge 0 \text{ and } \sum x_i = 1\}.$
$\cs.X$ has the usual face and degeneracy operations, and it is well known that its geometric realization $|\cs.X|$ has
the weak homotopy type of $X$.   See \cite{may} for details.  

Number the vertices of $\Delta^p$ $\{0, 1, \cdots, p\}$.  Then for a non-decreasing sequence of integers $0 \le i_0 \le
\ldots
\le i_r \le p$, define a   map of simplices    
$$(i_0, i_1, \ldots, i_r): \Delta^r \to \Delta^p,$$
by requiring that $(i_0, i_1, \ldots, i_r)$ be a linear map that  sends the vertex $k$ of $\Delta^r$ to the vertex $i_k$ of $\Delta^p$.
Given a singular $p$ - simplex $\sigma : \Delta^p \to X$, the composition of $\sigma$ with $(i_0, i_1, \ldots, i_r)$ defines an $r$ - simplex
$$
\sigma (i_0, i_1, \ldots, i_r) :  \Delta^r \to \Delta^p \to X.
$$

Now let $\begin{CD} F \hookrightarrow E @>\pi>> B \end{CD}$ be a fibration. There is a filtration of the simplicial set
$\cs.E$ defined as follows.

\begin{definition}  Let $F_p(\cs.E) \subset \cs.E$ be the subsimplicial set whose $r$ simplices are given by
 \begin{multline}
F_p (\cs_r(E)) =  \{T:\Delta^r \to E \, : \,  \pi  \circ T = \sigma (i_0, \ldots, i_r),  \,  \text{for some} \, 
\sigma \in \cs_q(B), \, q \leq p, \\
 \text{and some
sequence} \,  0\le i_0 \le \ldots \le i_r \le q\}. \notag
\end{multline}
\end{definition}
\med
\noindent Given a  simplicial set $Y.$, let $C_*(Y.)$ be the associated simplicial chain complex, whose $q$ - chains $C_q (Y.)$  are the free abelian group on
the
$q$ - simplices $Y_q$, and whose boundary homomorphisms are given by the alternating sum of the face maps.  Again, see \cite{may} for details. 
In particular for a space $X$, we have $C_*(\cs.X)$ is the singular chain complex which we previously denoted simply by $C_*(X)$.   The following is
is verified in \cite{MC}.

\begin{proposition}   Let $F \hk E \to B$ be a fibration as above.   Consider the   filtration of chain complexes,
$$
\{0\} \hk \cdots \hk F_{p-1}(C_*(E)) \hk F_p(C_*(E)) \hk \cdots \hk C_*(E)
$$
defined by $F_p(C_*(E)) = C_*(F_p(\cs.E))$ where $F_p(\cs.E)$ is the $p^{th}$ filtration of the singular simplicial set  defined above. 
Then this filtration induces the Serre spectral sequence converging to $H_*(E)$. 
\end{proposition}

\noindent We will study this spectral sequence in the examples of the fibrations
\begin{align}
ev : & LM \to M, \notag \\
\evi : & LM \times_M LM \to M, \notag \\
ev : & D_LM \to D_TM, \ \text{and} \notag \\
ev : & S_LM \to S_TM \notag
\end{align} described in the last
section.  In particular notice that by taking the filtration of pairs 
$$F_p(C_*(D_LM, S_LM)) = F_p(C_*(D_LM))/F_p(C_*(S_DM))$$
we get a spectral
sequence (the relative Serre spectral sequence)  converging to $H_*(D_LM, S_LM) = {\tilde H}_*(LM^{TM})$ and whose $E_2$ term is
$E^2_{p,q} = H_p(D_TM, S_TM; H_q(\Omega M))= {\tilde H}_p(M^{TM}; H_q(\Omega M)). $  Using the fibration $\evi : LM \times_M LM \to M$
there is a similar  filtration $F_p(C_*(\evi^*(D_TM), \evi^*(S_TM)))$ which yields a relative Serre spectral sequence converging to $\tilde H_*((LM
\times_M LM)^{TM})$ and whose $E_2$ term is ${\tilde H}_*(M^{TM}; H_*(\Omega M \times \Omega M))$. 

The following is an immediate observation based on the chain descriptions in the last section.

\med
\begin{proposition} \label{chain}The chain maps
\begin{align}
\times:\quad & C_*(LM) \otimes C_*(LM) \to C_*(LM \times LM) \quad  (\ref{c1}), \notag \\
\tilde{\tau}_\sharp:\quad & C_*(LM \times LM) \to C_*(\lmm)^{TM} \quad (\ref{c2}), \quad \text{and}\notag\\
\tilde{\gamma}_\sharp:\quad & C_*(\lmm)^{TM} \to C_*(LM^{TM}) \quad (\ref{c3}) \notag  
 \end{align} described in the last section  all preserve the above filtrations:
\begin{align}
\times:\quad & F_p(C_*(LM)) \otimes F_q(C_*(LM)) \to F_{p+q} (C_*(LM \times LM)),   \notag \\
\tilde{\tau}_\sharp:\quad & F_m(C_*(LM \times LM)) \to F_m (C_*(\evi^*(D_TM), \evi^*(S_TM))),\notag\\
\tilde{\gamma}_\sharp:\quad &  F_m (C_*(\evi^*(D_TM), \evi^*(S_TM)))\to F_m(C_*(D_LM, S_LM)), \notag  
 \end{align}
and therefore induce maps of the associated Serre spectral sequences.
\end{proposition}
\med
\noindent The following is a bit more delicate.
\med
\begin{theorem}\label{shift} The chain map
$$\tilde{u}_\sharp:   C_*(LM^{TM}) \to C_{*-d}(LM) \ (\ref{c4})$$
induces a map of filtered chain complexes  that lowers  the filtration by $d$,
$$
\tilde{u}_\sharp:  F_p(C_*(D_LM, S_LM))  \to F_{p-d} (C_{*}(LM)).
$$
and therefore induces a map of the associated Serre spectral sequences that  shifts grading,
$$
\tilde{u}_* : E^r_{p,q}(D_LM, S_LM) \to E^r_{p-d, q}(LM).
$$
In particular on the $E_2$ - level $\tilde u_*$ is the Thom isomorphism,
$$
H_p(M^{TM}; H_q(\Omega M)) \to H_{p-d}(M; H_q(\Omega M)).
$$
\end{theorem}

\med
\begin{proof}  
By the definition of the chain map $\tilde u_\sharp$, to prove this theorem it suffices to show that taking the cap product with the Thom class
$\ttt\in C^d(D_LM, S_LM)$ induces a map of filtered chain complexes that lowers the filtration by $d$, 
$$
\begin{CD}
F_p (C_*(D_LM, S_LM)) @>\cap \ttt >> F_{p-d} (C_*(D_LM)).
\end{CD}
$$

To verify this, recall that the cap product has the following chain level description.  Consider the operations on singular $n$ - simplices,
$$
\rfloor_p \quad \text{and} \quad {}_p\lfloor : \cs_n(X) \to \cs_p(X)
$$
defined by 
$$
\sigma\rfloor_p = (d_{p+1})^{n-p}(\sigma), \qquad {}_p\lfloor\sigma = (d_0)^{n-p}(\sigma).
$$
That is, $\sigma\rfloor_p $ is the restriction of $\sigma$ to the ``front"  $p$ - face, and ${}_p\lfloor\sigma $ is the restriction of $\sigma$ to the ``back"  $p$ -face.

Now we can choose our cochain $\ttt \in C^*(D_LM, S_LM )$  to represent the Thom class  so that  
 if  we view it  as an element in  
$\text{Hom}(C_d(D_LM);  \bz)$ it satisfies  
\begin{align}
& C_d(S_LM)  \subset \text{Ker}(  \ttt), \label{cap1} \\
& \text{Im}(s_j:  \cs_{d-1}(D_LM) \to  \cs_d(D_LM) \subset C_d(D_LM)) \subset \text{Ker} ( \ttt),
\label{cap2}
\end{align}
That is, $\ttt$ vanishes on chains on $S_LM$ and degenerate chains on $D_LM$.
Let $\sigma \in \cs_p (D_LM)$ for some $p \ge d$, then the cap product can be described as
\begin{equation}\label{ttt}
\sigma \cap \ttt = (-1)^{d(p-d)}\ttt({}_d\lfloor\sigma)  \cdot \sigma\rfloor_{p-d}.
\end{equation}
This gives a well defined chain map  $\cap \ttt: C_p(D_LM, S_LM) \to
C_{p-d}(D_LM)$ that represents capping with the Thom class in cohomology.  

The map $\cap \ttt$ on $C_*(D_LM, S_LM)$ is completely determined by its  composition with the projection
 $C_*(D_LM) \to C_*(D_LM,S_LM)$.  So to prove the lemma  it now   suffices to prove that $\cap \ttt$ maps $F_p(C_*(D_LM))$ to
$F_{p-d}(C_*(D_LM))$.

Let $T\in F_p(\cs_r(D_LM))$ be a singular $r$ - simplex in filtration $p$. Then, by definition, 
$$ev(T) = \sigma (i_0, i_1, \cdots i_r)$$ for some $q$-simplex $\sigma \in \cs_q(D_TM)$, and some
sequence $i_0 \leq \cdots  \leq i_r \leq q \le p$.   By the above formula for the cap product, we then have
$$
ev(T\cap \ttt) =  ev(T)\cap t = \pm t(\sigma (i_{r-d}, \ldots, i_r))\sigma(i_0, \ldots, i_{r-d}),
$$ where, as above, $t \in C^d(D_TM, S_TM)$ represents the Thom class of  the tangent bundle $TM \to M$.  By  (\ref{cap2}) $t$
vanishes on degenerate simplices.  Therefore this expression can only be nonzero if $i_{r-d} < \ldots < i_r \leq q$, and therefore $r-d \leq q-d \le p-d$.
Hence $T\cap \ttt \in F_{p-d} (C_*(D_LM))$ as claimed. \end{proof}

\med
\noindent Notice that if we  compose the chain maps in proposition \ref{chain} and theorem \ref{shift}, we have a map of filtered chain
complexes 
$$\mu : F_p(C_*(LM)) \otimes F_q(C_*(LM)) \la F_{p+q-d}(C_*(LM))$$
which,  by proposition \ref{loopprod} induces the loop product in homology.  Therefore $\mu$ induces a map of spectral sequences
\begin{equation}\label{mu}
\mu : E^r_{p, s}(LM) \otimes E^r_{q,t}(LM)  \to  E^r_{p+q-d, s+t}(LM).
\end{equation}
For simply connected $M$, on the $E_2$ - level, $\mu$ defines a map
$$
\mu : H_p(M; H_s(\Omega M)) \otimes H_q(M; H_t(\Omega M)) \la H_{p+q-d}(M; H_{s+t}(\Omega M))
$$ which we claim is given up to sign, by the intersection product with coefficients on the Pontryagin ring $H_*(\Omega M)$. More explicitly,
$$
\mu((a \otimes g) \otimes (b \otimes h)) = \pm (a \cdot b)  \otimes (gh)
$$
where $a \in H_p(M)$,  $g \in H_s(\Omega M)$,  $b \in H_q(M)$, and $h \in H_t(\Omega M)$, and where    
$a \cdot b$ is the intersection product and $gh$ is the Pontryagin product.

To see this, notice that the composition of chain maps used to define $\mu$
is given by a composition of chain maps of fibrations  (and pairs of fibrations).  On the   base space level this is given by
the composition of maps described in proposition \ref{intersection} realizing the intersection product.  On the fiber level, the
fact that this chain map induces the Pontrjagin product comes from the fact that  map $\gamma : LM \times_M LM \to LM$ 
as defined in (\ref{gamma}) is a map of fibrations

$$
\begin{CD}
\Omega M \times \Omega M @>\rho >>  \Omega M \\
@VVV   @VVV \\
LM \times_M LM @>\gamma >>  LM \\
@V \evi VV  @VVev V \\
M  @>> = > M
\end{CD}
$$ where $\rho$ is the Pontryagin product. 

Thus $\mu$ defines a multiplicative structure on the Serre spectral sequence for the fibration $\Omega M \to LM \to M$  which converges to the loop homology algebra structure
on $H_*(LM)$,  and on the $E^2$ -level is given  (up to sign) by the intersection product on $M$  with coefficients in the Pontryagin ring $H_*(\Omega M)$.  
 However the grading is shifted in a
way that is confusing for calculational purposes.  To remedy this, define a second quadrant spectral $\{ E^r_{s.t}
(\bh_*(LM));  \, d_r : E^r_{s,t} \to E^r_{s-r, t+r-1}\}$ with
$s
\leq 0$,
$t
\geq 0$, which converges to the loop homology  $\bh_{s+t} (LM)$   by
regrading the Serre spectral spectral sequence  in the following way. 
\begin{equation}\label{dshift}
 E^r_{s.t} (\bh_*(LM))=  E^r_{s+d, t}(LM),
\end{equation}
where $d$ is the dimension of $M$, and the right hand side is the Serre spectral sequence for the fibration $\Omega M \to LM \to M$ we have been
considering.   Notice that $E^r_{s.t} (\bh_*(LM))$ can only be nonzero for $-d \leq s \leq 0$. Moreover with the new indexing the spectral sequence converges to
the loop homology in a grading preserving way,  $E^r_{s.t} (\bh_*(LM)) \rightrightarrows \bh_{s+t}(LM)$. We  also see that
with this new indexing,  the  loop multiplication in the spectral sequence (\ref{mu}) preserves the bigrading,  
$$
\mu : E^r_{s,t}(\bh_*(LM)) \otimes E^r_{p,q}((\bh_*(LM)) \la E^r_{s+p, t+q}((\bh_*(LM)).
$$
Finally notice that the $E^2$ - term is given by
$$
 E^2_{s,t}(\bh_*(LM)) = H_{s+d}(M; H_t(\Omega M))
$$
for $-d \leq s \leq 0$.  By applying Poincare duality we have
$$
E^2_{s,t}(\bh_*(LM)) = H^{-s}(M; H_t(\Omega M)).
$$
Since under Poincare duality the intersection pairing in homology  coincides with the cup product in cohomology,
on the $E^2$ - level the multiplication $\mu$  is given (up to sign) by cup product  with coefficients in the Pontryagin ring $H_*(\Omega M)$.
$$
\mu :  H^{-s}(M; H_t(\Omega M)) \otimes H^{-p}(M; H_q(\Omega M)) \to H^{-(s+p)}(M; H_{t+q}(\Omega M))
$$
for $-d \leq s, p \leq 0$ and $t, q \geq 0$.
This completes the proof of theorem \ref{ss}.
\med
\section{The Loop Product on $S^n$ and $\bcp^n$}

\med
\noindent Our goal in this section is to use the loop homology spectral sequence constructed in the last section to perform the
calculations described in theorems \ref{spheres} and \ref{projspace}.

 We first prove theorem \ref{spheres} by  calculating  the ring structure of $\bh_*(LS^n)$. If $n=1$,  the fibration
$\Omega S^1 \to LS^1 \to S^1$ is trivial.  Since the components of the based loop space $\Omega S^1$ are all contractible,
and $\pi_0(\Omega S^1) \cong \bz$,  we have a homotopy equivalence
$$LS^1 \simeq S^1 \times \Omega S^1 \simeq S^1 \times \bz.$$
Similarly, there is a   homotopy equivalence $LS^1 \times_{S^1} LS^1 \simeq S^1 \times \bz \times \bz$.  With respect to these equivalences, it is clear
that the map
$\gamma : LS^1 \times_{S^1} LS^1  \to LS^1$ is  given by
\begin{align}
S^1 \times \bz \times \bz &\la S^1 \times \bz \notag \\
(x,m,n) &\la (x, m+n).
\end{align}
It is also clear that with respect to these equivalences, the Thom - Pontryagin map
$\tilde \tau : LS^1 \times LS^1 \to  (LS^1 \times_{S^1} LS^1)^{TS^1}$ is given by
$$
\tau \times 1 : S^1 \times S^1 \times \bz \times \bz  \la (S^1)^{TS^1} \wedge (\bz \times \bz)_+
$$
where $\tau : S^1 \times S^1 \to (S^1)^{TS^1}$ is the Thom - Pontryagin construction for the diagonal map $\Delta : S^1 \to S^1 \times S^1$.  
   Thus by (\ref{prod}) the loop homology algebra structure on $\bh_*(LS^1)$  is, with respect to the equivalence $LS^1 \simeq S^1 \times \bz$,  given
by the tensor product of the intersection ring structure on $H_*(S^1)$ with the group algebra structure, $H_0(\bz) \cong \bz[t,t^{-1}]$.  Using Poincare
duality, we then have an algebra isomorphism 
\begin{equation}\label{s1}
\bh_{-*}(LS^1) \cong H^*(S^1) \otimes H_0 (\bz) \cong   \Lambda [a]  \otimes \bz[t, t^{-1}]
\end{equation}
where $a \in \bh_{-1}(LS^1)$ corresponds to the generator in $H^0 (S^1)$.  

We now proceed with a calculation of $\bh_*(LS^n)$ for $n > 1$.   Consider the loop homology spectral sequence  in this case.
For dimension reasons, the only nontrivial differentials occur at  the $E^n$ level. Recall that there is an isomorphism of algebras,
$H_*(\Omega S^n)
\cong
\bz[x]$,  where $x$ has degree
$n-1$. It then follows that
$$E^2_{-p, q}(\bh_*(LS^n)) \cong \cdots \cong E^n_{-p, q} \cong H^p(S^n) \otimes H_q(\Omega S^n).$$

The differentials $d^n$ can be computed using the results in \cite{B} and \cite{S}. An exposition of this calculation (for the
Serre spectral sequence of the fibration $\Omega S^n \to LS^n \to S^n$) is given in
\cite{M}.  Inputting the change of grading used to define the loop homology spectral sequence, we have the following
picture of the differentials: 

\begin{center}
\begin{picture}(360, 230)(-20, -25)
\put(-20, 190){\em $E^n$ terms:}
\put(105, 0){\vector(-1, 0){65}}
\put(105, 0){\vector(0, 1){185}}
\put(-2, 0){$H^{-*}(S^n)$}
\put(85, 190){$H_*(\Omega S^n)$}
\multiput(65, 0)(0, 30){5}{\circle*{3}}
\multiput(105, 0)(0, 30){5}{\circle*{3}}
\multiput(105, 0)(0, 30){4}{\vector(-4, 3){40}}
\multiput(80, 20)(0, 30){4}{$\times 0$}
\multiput(80, 140)(0, 15){3}{$\cdot$}
\multiput(108, 3)(0, 30){5}{$\Z$}
\multiput(57, 3)(0, 30){5}{$\Z$}
\put(10, 160){$n = \text{odd}$}
\put(62, -9){$-n$}
\put(127, 27){$n-1$}
\put(123, 30){\vector(-1, 0){14}}

\put(280, 0){\vector(-1, 0){65}}
\put(280, 0){\vector(0, 1){185}}
\put(173,0){$H^{-*}(S^n)$}
\put(260, 190){$H_*(\Omega S^n)$}
\multiput(240, 0)(0, 30){5}{\circle*{3}}
\multiput(280, 0)(0, 30){5}{\circle*{3}}
\multiput(280, 0)(0, 30){4}{\vector(-4, 3){40}}
\multiput(255, 50)(0, 60){2}{$\times 2$}
\multiput(255, 20)(0, 60){2}{$\times 0$}
\multiput(255, 140)(0, 15){3}{$\cdot$}
\multiput(283, 3)(0, 30){5}{$\Z$}
\multiput(232, 3)(0, 30){5}{$\Z$}
\put(185, 160){$n = \text{even}$}
\put(237, -9){$-n$}
\put(302, 27){$n-1$}
\put(298, 30){\vector(-1, 0){14}}

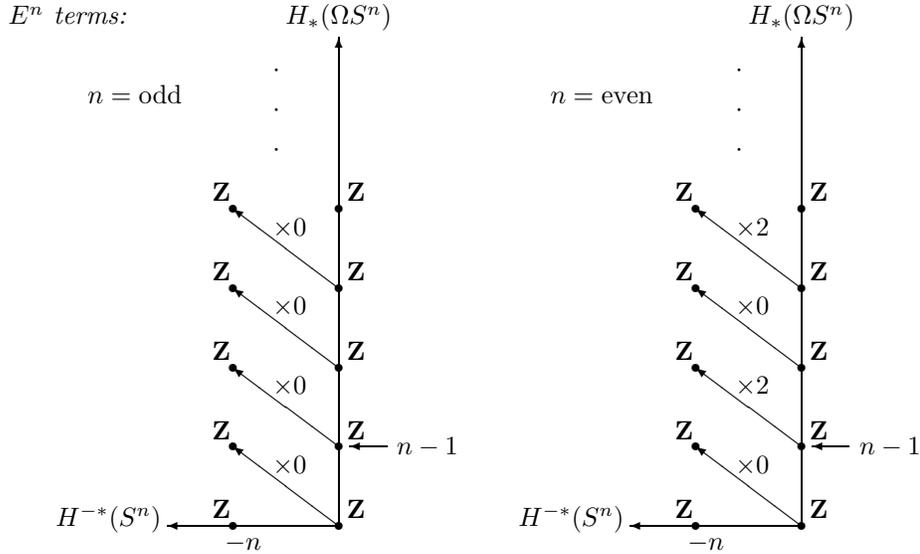
\put(80, -22){The Spectral Sequences for $LS^n$}
\end{picture}
\end{center}

Denote by $\iota$, $\sigma$ the generators of $H^n(S^n)$ and $H^0(S^n)$, respectively. Let $1_\Omega$ be the unit of 
$H_*(\Omega S^n)$. In both spectral sequences $\sigma \otimes 1_\Omega$ is an infinite cycle and therefore represents a
class in
$\bh_0 (LS^n)$. This class is the unit in the algebra $\bh_*(LS^n)$, which we denote by $1$. Similarly $a = \iota \otimes
1_\Omega$ represents a class in $\bh_{-n}(LS^n)$,  and for dimension reasons $a^2$ must vanish in $\bh_*(LS^n)$.

Let $u = \sigma \otimes x \in E^n_{0, n-1}$.  In the case when $n$ is odd,  the spectral sequence collapses. 
Since $n > 1$, for dimensional reasons there can be no extension problems, and so we have an isomorphism of algebras
$$
 \bh_*(LS^n) \cong E^2_{*, *} \cong H^{-*}(S^n) \otimes H_*(\Omega S^n) \cong \Lambda[a] \otimes \bz[u].   $$

If $n$ is even,    the $u^{2k}$ terms   survive to $E^\infty$, and   
$$
d_n(u^{2k+1}) =   2au^{2k+2}
$$
for all nonnegative integers $k$.    Let $v = u^2 = \sigma \otimes x^2$.  Being an infinite cycle it represents a class in
$\bh_{2(n-1)}(LS^n)$.  Then
$v^k$ represents a generator of a subgroup  of $\bh_{2k(n-1)}(LS^n)$ represented by classes in $E^\infty_{0, 2k(n-1)}$. 
 Notice that the ideal
$(2av)$ vanishes in
$E^\infty$.   Thus $\bz[a, v]/(a^2, 2av)$ is a subalgebra in $ E^\infty_{*,*}$.
Let $b  = \iota \otimes x \in E^\infty_{-n, n-1}$.   Then $E^\infty_{*,*} $ is generated by
$\bz[a, v]/(a^2, 2av)$ and
$b$. Again for dimension reasons $b^2 = ab = 0 \in E^\infty$,  and so
$$E^\infty_{*,*} \cong (\Lambda[b] \otimes \bz[a, v])/(a^2, ab, 2av), \qquad \text{if $n$ is even.}$$
When $n > 2$, for dimensional reasons there can be no extension issues, so
$\bh_*(LS^n) \cong (\Lambda[b] \otimes \bz[a, v])/(a^2, ab, 2av)$ for $n$ even and $n > 2$.   For $n=2$ we consider the potential extension problem.
For filtration reasons, there are unique classes in $\bh_*(LS^2)$ represented by $a$ and $b$ in $E^\infty_{*,*}$.  The ambiguity in the choice
of class represented by $v \in E^\infty_{0,2}$ lies in $E^\infty_{-2,4} \cong \bz$ generated by $av$.   Since $a^2=0$ in $\bh_*(LS^2)$, any choice
$\tilde v \in \bh_2(LS^2)$ will satisfy $2a\tilde v = 0$. Thus any choice of $\tilde v$ together with $a$ and $b$  will generate the same algebra, namely
 $\bh_*(LS^2) \cong (\Lambda[b] \otimes \bz[a, \tilde v])/(a^2, ab, 2a\tilde v)$. 

This completes the proof of theorem \ref{spheres}. 

\med
\noindent We now proceed with the proof of theorem \ref{projspace} by calculating
  the ring structure of $\bh_*(\bcp^n)$. Notice that if $n=1$, $\bh_* (L\bcp^n) \cong \bh_* (LS^2)$, and this
case was  already discussed   above.  So for what follows we   assume $n > 1$.  The $E^2$ -term of the loop
homology spectral sequence is
$H^*(\bcp^n; H_*(\Omega \bcp^n))$.  The cohomology ring
$H^*(\bcp^n) \cong \bz[c_1]/(c_1^{n+1})$ is generated by the first Chern class $c_1$.  We now recall the Pontryagin ring structure of
  $H_* (\Omega \bcp^n)$. 

Consider the homotopy fibration 
$$
\begin{CD}
\Omega S^{2n+1} @>\Omega\eta >>  \Omega \bcp^n  \to \Omega \bcp^\infty \simeq S^1,
\end{CD}
$$
where $\eta : S^{2n+1} \to \bcp^n$ is the Hopf map.  Since this is a fibration of loop spaces that has a section (because $\pi_1(\Omega \bcp^n) \cong
\pi_1(\Omega \bcp^\infty))$ we can conclude the following.
\begin{enumerate}
\item The fibration is homotopically trivial, $\Omega \bcp^n \simeq \Omega S^{2n+1} \times S^1$.   
\item The Serre spectral sequence of this fibration is a spectral sequence of algebras.
\end{enumerate}
Therefore we have that in the Serre spectral sequence,
\begin{align}E^\infty_{*,*} &\cong H_*(S^1) \otimes H_*(\Omega S^{2n+1}) \notag \\
&\cong \Lambda[t] \otimes \bz[x]
\end{align}
where $t \in H_1(S^1)$ and $x \in H_{2n}(\Omega S^{2n+1})$ are generators.  This isomorphism is  one of algebras.  Indeed it is clear that there are no extension
issues in this spectral sequence so that 
$$
H_*(\Omega\bcp^n) \cong H_*(S^1) \otimes H_*(\Omega S^{2n+1}) \cong \Lambda[t] \otimes \bz[x].
$$

Now consider the loop homology spectral sequence for $L\bcp^n$. We therefore have
\begin{align}
E^2_{-*,*}(\bh_*(L\bcp^n))  &\cong  H^*(\bcp^n; H_*(\Omega \bcp^n)) \notag \\
&\cong \bz[c_1]/(c_1^{n+1})  \otimes \Lambda[t] \otimes \bz[x] \notag
\end{align}
as algebras. 

It was computed in \cite{Z} that
$$H_k(L\bcp^n) = \begin{cases}
\Z \qquad\qquad & \text{if $k=0, 1, \ldots, k \ne 2mn$, $m \ge 1$}, \\
\Z \oplus \Z_{n+1} & \text{if $k=2mn$, $m \ge 1$}.\end{cases}
$$
This calculation implies the following pattern of differentials in the loop homology spectral sequence:
\begin{center}
\begin{picture}(360, 240)(-60, -30)
\put(-20, 190){\em $E^{2n}$ terms:}
\put(205, 0){\vector(-1, 0){135}}
\put(205, 0){\vector(0, 1){185}}
\put(20, 0){$H^{-*}(\bcp^n)$}
\put(180, 190){$H_*(\Omega \bcp^n)$}
\multiput(205, 0)(0, 65){3}{\circle*{3}}
\multiput(205, 5)(0, 65){3}{\circle*{3}}
\multiput(185, 0)(0, 65){3}{\circle*{3}}
\multiput(185, 5)(0, 65){3}{\circle*{3}}
\multiput(165, 0)(0, 65){3}{\circle*{3}}
\multiput(165, 5)(0, 65){3}{\circle*{3}}
\multiput(145, 0)(0, 65){3}{\circle*{3}}
\multiput(145, 5)(0, 65){3}{\circle*{3}}
\multiput(120, 0)(0, 65){3}{$\cdots$}
\multiput(105, 0)(0, 65){3}{\circle*{3}}
\multiput(105, 5)(0, 65){3}{\circle*{3}}
\multiput(85, 0)(0, 65){3}{\circle*{3}}
\multiput(85, 5)(0, 65){3}{\circle*{3}}
\multiput(205, 5)(0, 65){2}{\vector(-2, 1){120}}
\multiput(95, 145)(0, 10){3}{$\cdot$}
\multiput(125, 145)(0, 10){3}{$\cdot$}
\multiput(155, 145)(0, 10){3}{$\cdot$}
\multiput(175, 145)(0, 10){3}{$\cdot$}
\multiput(195, 145)(0, 10){3}{$\cdot$}
\put(40, 160){$n > 1$}
\put(79, -10){$-2n$}
\put(203, -10){$0$}
\put(183, -10){$-2$}
\put(163, -10){$-4$}
\put(143, -10){$-6$}
\multiput(144, 40)(0, 65){2}{$\times (n+1)$}
\put(209, -2){${}_0$}
\put(209, 4){${}_1$}
\put(209, 63){${}_{2n}$}
\put(209, 69){${}_{2n+1}$}
\put(209, 128){${}_{4n}$}
\put(209, 134){${}_{4n+1}$}
\put(45, -25){The  loop homology spectral sequence for $L\bcp^n$}
\end{picture}
\end{center}

That is,  the only nonzero differentials are $d_{2n}$ and
$$d_{2n}(y) =   (n+1) c^{n} u,$$
where $c = c_1 \otimes 1$, $u = \sigma \otimes x$, and $y = \sigma \otimes t$. 
Here $\sigma$ is a  generator of $H^0 (\bcp^n) \cong \bz$.  For
dimension reasons
$d_{2n} (u) = 0$,  and since $d_{2n}$ is a derivation,
$$d_{2n}(yu^k) = d_{2n}(y) u^k = (n+1) c^n u^{k+1},$$
and these are all the nonzero differentials.

The spectral sequence collapses beyond $E^{2n}$ level,  therefore $c$ and $u$ represent homology classes in $\bh_{-2}(L\bcp^n)$ and $\bh_{2n}(L\bcp^n)$
respectively. Notice also that the ideal
$(c^{n+1}, (n+1)c^nu)$ vanishes in $E^\infty_{*, *}$.   Therefore we have  a subalgebra
$$\Z[c, u]/(c^{n+1}, (n+1)c^n u) \subset E^\infty_{*,*}(\bh_* (L\bcp^n)).$$
Let $w = yc \in E^2_{*,*}.$  $w$ is an infinite cycle in this spectral sequence and represents a class in $\bh_{-1}(\bcp^n)$.  Notice  that $w^2 = 0$ in
$E^\infty_{-4, 2}$. Similarly $wc^n$ also  vanishes in $E^\infty_{*,*}$.  Thus the $E^\infty$ - term of the loop homology spectral sequence can be written as
follows:
$$E^\infty_{*,*}(\bh_* (L\bcp^n) )\cong (\Lambda[w] \otimes \Z[c, u])/(c^{n+1}, (n+1)c^n u, wc^n).$$  The extension issues are handled just as they were for
$\bh_*(LS^2) = \bh_*(L\bcp^1)$ above.  We then conclude that
$$ \bh_* (L\bcp^n) \cong (\Lambda[w] \otimes \Z[c, u])/(c^{n+1}, (n+1)c^n u, wc^n)$$  as algebras.  This completes the proof of theorem \ref{projspace}.

\med

\end{document}